\documentclass[12pt]{amsart}
\usepackage{latexsym}
\usepackage{amsmath}
\usepackage{amsfonts}
\usepackage{amssymb}
\usepackage{amsthm}

\theoremstyle{plain}
\newtheorem{theorem}{Theorem}[section]
\newtheorem{lemma}[theorem]{Lemma}
\newtheorem{proposition}[theorem]{Proposition}

\newtheorem*{Main Theorem}{Main Theorem}

\newtheorem{corollary}[theorem]{Corollary}

\theoremstyle{definition}
\newtheorem{definition}[theorem]{Definition}

 \theoremstyle{remark}
\newtheorem*{remark}{Remark}
\newtheorem*{remarks}{Remarks}

\newcommand{\N}{\mathbb{N}}

\newcommand{\R}{\mathbb{R}}
\newcommand{\T}{\mathbb{T}}
\newcommand{\Z}{\mathbb{Z}}
\newcommand{\A}{\mathcal{A}}
\newcommand{\B}{\mathcal{B}}
\newcommand{\C}{\mathcal{C}}
\newcommand{\D}{\mathcal{D}}
\newcommand{\E}{\mathcal{E}}
\newcommand{\F}{\mathcal{F}}

\newcommand{\I}{\mathcal{I}}
\newcommand{\M}{\mathcal{M}}

\newcommand{\bX}{{\bf X}}
\newcommand{\bY}{{\bf Y}}
\newcommand{\bD}{{\bf D}}
\newcommand{\bS}{{\bf S}}
\newcommand{\bT}{{\bf T}}

\newcommand{\gG}{\Gamma}

\newcommand{\bbE}{\mathbb{E}}

\renewcommand{\O}{\mathcal{O}}
\renewcommand{\P}{\mathcal{P}}

\newcommand{\e}{\varepsilon}
\newcommand{\V}{\bigvee}

\renewcommand{\)}{\right)}
\renewcommand{\(}{\left(}

\def \colon{{:}\;}
\begin{document}

\title{The Structure of Strongly Stationary Systems}
\author{NIKOS FRANTZIKINAKIS}

\address{ Department of Mathematics, McAllister Building,
Pennsylvania State University,
 University Park, PA 16802} \email{nikos@math.psu.edu}
\pagestyle{headings}

\begin{abstract}
 Motivated by a problem in ergodic Ramsey theory, Furstenberg and Katz-\\
 nelson
 introduced the notion of strong stationarity,
showing that certain recurrence properties hold for arbitrary
measure preserving systems if they are valid for strongly
stationary ones. We construct some new examples and prove a
structure theorem for strongly stationary systems. The building
blocks are Bernoulli systems and rotations on nilmanifolds.
\end{abstract}

\maketitle

\section{Introduction}
\subsection{Historical background.}
In 1975 Szemer\'edi proved the following long standing conjecture
of Erd\"os and Tur\`an:
\begin{theorem}[{\bf Szemer\'edi}]
Let $\Lambda$ be a subset of the integers with positive upper
density. Then $\Lambda$ contains arbitrarily long arithmetic
progressions.
\end{theorem}
Szemer\'edi's proof was combinatorial in nature and intricate. In
$1977$ Furstenberg  (\cite{Fu1})  gave an entirely different proof
using ergodic theory. He showed that Szemer\'edi's theorem is
equivalent to a statement about multiple recurrence of measure
preserving systems and then proved the ergodic version:

\begin{theorem}[{\bf Furstenberg}]\label{T:multrec}
 Let $(X,\B,\mu,T)$ be a finite measure preserving system
and $A\in \B$ be a set
 with positive measure. Then for every $k\in \N$, there exists $n\in\N$
 such that
 $$
 \mu(A\cap T^{-n}A\cap\cdots \cap T^{-nk}A)>0.
 $$
\end{theorem}

Furstenberg's proof  launched the field of ergodic Ramsey theory,
where problems in combinatorics are translated to recurrence
related statements of measure preserving systems and then proved
using ergodic theory. Using this approach Furstenberg and
Katznelson, and more recently Bergelson and Leibman (among others)
established several ergodic theoretic results whose combinatorial
implications are not currently attainable by any other methods.
This includes a multidimensional and a polynomial extension of
Szemer\'edi's theorem (\cite{FK1}, \cite{BL}), and the density
version of Hales-Jewett theorem (\cite{FK2}), the "master" theorem
that contains several  major results in the field as special
cases.

 The notion of strong stationarity (Definition \ref{D:sst}) was
introduced in the paper of Furstenberg and Katznelson (\cite{FK2})
in proving the density version of  Hales-Jewett theorem. An
important result  established in the same paper is that an
arbitrary stationary process "majorizes"  a strongly stationary
one (Section \ref{SS:majorize}). From this it follows  that
several recurrence properties are valid for arbitrary measure
preserving systems if they are valid for strongly stationary ones.
In particular, it turns out to be sufficient to verify Theorem
\ref{T:multrec} for the specific case of a strongly stationary
system. This motivated the problem of determining the structure of
strongly stationary systems.

In \cite{Je2} Jenvey proved that every  ergodic strongly
stationary system is necessarily Bernoulli. Unfortunately, not
every ergodic system majorizes an ergodic strongly stationary one,
and the ergodic components of a strongly stationary system are not
necessarily strongly stationary. In fact  nonergodic strongly
stationary systems can have completely different structure than
the ergodic ones; there exist several distal examples that one has
to identify.

We give a structure theorem for strongly stationary systems
(Theorems \ref{T:components} and \ref{T:main}):
\begin{Main Theorem}
(i) Almost every  ergodic component of a strongly stationary
system is isomorphic to the direct product of a Bernoulli system
and a totally ergodic pro-nilsystem (defined
in Section \ref{SS:nilsystems}).\\
(ii) An extremal strongly stationary system (Definition
\ref{D:extremal})  is isomorphic to the direct product of a
Bernoulli system and a strongly stationary system associated to
some pro-nilmanifold (defined in Section \ref{SS:last}).
\end{Main Theorem}

Moreover, we construct new examples of strongly stationary systems
(Section \ref{SS:last}, examples (iv) and (v)).

\subsection{Format of the paper.}

After reviewing some preliminary notions and results in Section 2
we define strong stationarity in Section 3. We prove that an
arbitrary  stationary process majorizes a strongly stationary one
and  give the basic examples of strongly stationary systems.

The general strongly stationary system is an integral of extremal
ones. In Section 4 we give necessary and sufficient conditions for
extremality and prove a homogeneity property for extremal strongly
stationary systems.

 In Section 5 we prove that almost every ergodic component of a strongly
stationary system is isomorphic to the direct product of a distal
system and a Bernoulli system. Moreover, we show that the distal
factor of a strongly stationary system is strongly stationary and
coincides with the characteristic factor of the system (defined in
Section \ref{SS:cfactor}). This reduces our problem to determining
the structure of distal strongly stationary systems.

Finally, in Section 6 we obtain a structure theorem for distal
strongly stationary systems. Using  results from Section 5, in
conjunction with a recent result of Host and Kra (\cite{HK3}), we
show that their ergodic components are pro-nilsystems. Moreover,
we construct new examples of distal strongly stationary systems
with nonaffine ergodic components. This new set of examples allows
us to give a complete classification.

\section{Preliminaries} To facilitate the reading, we establish
our notation and review some basic results that are used in the
sequel. We refer the reader to  \cite{Fu2}, \cite{Pe}, and
\cite{Wa} for more details.

\subsection{Measure preserving systems.} A \emph{measure preserving
system} (or just system) is a measure space $(X,\B,\mu)$ together
with a measurable measure preserving transformation $T$ on it.
Throughout the discussion we assume that all measure spaces are
Lebesgue. When there is no danger of confusion we use the bold
symbol ${\bf X}$ to denote the system $(X,\B,\mu,T)$. We also use
the bold symbol $\bT$ to denote the operator $\bT\colon
L^{\infty}({\bf X})\to L^{\infty}({\bf X})$ defined by $(\bT
f)(x)= f(T(x))$.

Let ${\bf X}=(X,\B,\mu,T)$ be  an invertible measure preserving
system and $\B_0$ be a sub-$\sigma$-algebra of $\B$. By $\V_{-n}^m
T^i \B_0$ we denote the $\sigma$-algebra spanned by sets of the
form $T^i B$ where $B\in \B_0$,  $-n\leq i \leq m$. We say that
$\B_0$ is $T$-\emph{generating} if $\V_{-\infty}^\infty T^i
\B_0=\B$ up to null sets. If $\F$ is an algebra of bounded
$\B$-measurable complex valued functions,  we denote by $\B(\F)$
the sub-$\sigma$-algebra of $\F$-measurable sets, that is, the
$\sigma$-algebra generated by sets of the form $f^{-1} (A)$, where
$f\in\F$ and $A\subset \mathbb{C}$ is open. We say that $\F$ is
$T$-\emph{generating} if $\B(\F)$ is $T$-generating.

\subsection{ Furstenberg's structure theorem.}
\label{SS:furstenberg} Let ${\bf X}$ be a measure preserving
system,  ${\bf Y}=(Y,\mathcal{E},\nu,R)$ be a factor of ${\bf X}$,
and $\mu=\int\mu_y \ d\nu(y)$ be the disintegration of $\mu$ over
${\bf Y}$. A function $f\in L^2({\bf X})$ is \emph{compact}
relative to the factor ${\bf Y}$ if for every $\e>0$ there exist
functions $g_1,\ldots,g_m\in L^2({\bf X})$ such that $\min_{1\leq
s\leq m}$ $\|\bT^if-g_s\|_{L^2(\mu_y)}<\e$ for every $i\in\N$ and
for $\nu$-a.e. $y\in Y$. An extension ${\bf X}$ of ${\bf Y}$ is
\emph{compact} if the set of compact functions relative to ${\bf
Y}$ is dense in $L^2({\bf X})$.

Starting with the factor of $T$-invariant functions ${\bf D}_0$,
we define inductively ${\bf D}_{i+1}$ to be the maximal compact
extension of ${\bf D}_i$. More precisely, we consider the
subalgebra  generated by the compact functions relative to ${\bf
D}_i$ and  we define ${\bf D}_{i+1}$ to be the factor determined
by it. We call the factor ${\bf D}_k$ the \emph{$k$-step distal
factor} of the system ${\bf X}$. The maximal factor that can be
exhausted by a transfinite number of compact extensions is called
the \emph{distal factor} and is denoted by ${\bf D}$.

 We say that ${\bf X}$ is a \emph{relatively ergodic} extension
of ${\bf Y}$, if every $T$-invariant function on  ${\bf X}$ is
$\mu$-a.e. a function on ${\bf Y}$ (that is $\E$-measurable).

Let ${\bf X}_1$, ${\bf X}_2$ be two extensions of ${\bf Y}$. The
fiber product space is defined as
$$
X_1\times_Y X_2=\left\{(x_1,x_2)\in X_1\!\times\! X_2\colon
\pi_1(x_1)=\pi_2(x_2) \right\},
$$
where $\pi_i\colon X_i\to Y$ the factor map. By \ $\B'$ we denote
the restriction of $\B_1\times\B_2$ on $ X_1\!\times_Y\! X_2$, and
by $\mu'$ the measure  defined by the disintegration (over ${\bf
Y}$) $\mu'_y=\mu^1_y\times\mu^2_y$, where $\mu_1=\int \mu^1_y \
d\nu(y)$ and $\mu_2=\int \mu^2_y \ d\nu(y)$. If ${\bf X}_1={\bf
X}_2={\bf X}$, and ${\bf Y}$ is the factor determined by the
action of $T$ on a sub-$\sigma$-algebra $\E$ of $\B$, then
$$
\int f(x_1) f(x_2) \ d\mu'=\int \bbE(f_1|\E) \bbE(f_2|\E) \ d\mu
$$
where $\bbE(f|\E)$ denotes the conditional expectation of $f$
given $\E$.
 We can check that ${\bf X}_1\!\times_{\bf Y}\! {\bf X}_2$$=(X_1\!\times_Y\!
X_2, \B',\mu',T_1\!\times\! T_2)$ is a measure preserving system
that extends ${\bf Y}$. We say that ${\bf X}$ is a \emph{relative
weak mixing} extension of ${\bf Y}$ if ${\bf X}\!\times_{\bf Y}\!
{\bf X}$ is a relatively ergodic extension of ${\bf Y}$.

Furstenberg (\cite{Fu2}) proved the following structure theorem:
\begin{theorem}[{\bf Furstenberg}]
Every  measure preserving system is  a relative weak mixing
extension of its distal factor.
\end{theorem}
We note that the distal factor is the smallest factor with respect
to which the system is relative weak mixing.
The  following result was proved in \cite{Fu1} for  ergodic
systems and is needed in the sequel. The argument given there also
works for nonergodic systems.
\begin{lemma}[{\bf Furstenberg}]\label{L:kronecker}
The invariant functions of  \ ${\bf X\times_Y X}$  belong to the
closed subspace spanned by functions of the form
$f_1(x_1)f_2(x_2)$ where $f_1,f_2$ are  functions compact relative
to ${\bf Y}$.
\end{lemma}

\subsection{Stationary processes.} A $\Lambda$-valued
\emph{stochastic process} is a sequence of measurable functions
(random variables) $\{f_i\}_{i\in \Z}$ defined on a probability
space $(X,\B,\mu)$ with values in a compact metric space $\Lambda$
(the state space).

The \emph{finite dimensional statistics} of a stationary process
$\{f_i\}_{i\in \Z}$ is the collection of all measurements $
\mu\big(\bigcap_{i=-k}^{k} \{f_i \in A_i\}\big)$, where $k\in \N$
and $A_i\subset\Lambda$ are open.

A  stochastic process is  \emph{stationary} if its finite
dimensional statistics are invariant under translations of the
time parameter, that is, $ \mu\big(\bigcap_{i=-k}^{k} \{f_{i+r}
\in A_i\}\big)= \mu\big(\bigcap_{i=-k}^{k} \{f_{i} \in A_i\}\big)$
for all $k,r\in \N$ and open sets $A_i\subset \Lambda$.

Two stationary processes $\{f_i\}_{i\in \Z}$ and $\{g_i\}_{i\in
\Z}$ (suppose that $\mu,\nu$ are the corresponding underlying
measures) are \emph{equivalent} if they have the same finite
dimensional statistics, that is, $ \mu\big(\bigcap_{i=-k}^{k}
\{f_{i} \in A_i\}\big) = \nu\big(\bigcap_{i=-k}^{k} \{g_{i} \in
A_i\}\big)$ for every $k \in \N$ and open sets $A_i \subset
\Lambda$.


\subsection{Sequence space representations.} Let $\Lambda$ be a
compact metric space. The \emph{sequence space} $\Lambda^\Z$
equipped with the product topology is again a compact metric
space. We denote by $x_i$ the $i$-th coordinate of a point $x\in
\Lambda^\Z$. The  Borel $\sigma$-algebra ${\A}^\Z$ is generated by
the finite dimensional rectangles $\bigcap_{-k}^k\{x\colon x_i \in
A_i \}$ where each $A_i\subset \Lambda$ is open. A probability
measure $\sigma$ defined on the completion of $\A^\Z$ (which we
denote again by $\A^\Z$) is \emph{stationary} if $
\sigma\big(\bigcap_{i=-k}^k\{x\colon x_i \in A_i \}\big)=
\sigma\big(\bigcap_{i=-k}^k\{x\colon x_{i+r} \in A_i \}\big)$ for
every $k \in \N$, $r \in \Z$ and open sets $A_i\subset\Lambda$.
Having fixed the space $\Lambda$ we denote by $\M$ the set of all
probability measures and by $\M_s$ the space of all stationary
measures on the sequence space $\Lambda^\Z$. Both $\M$ and $\M_s$
endowed with the weak-star ($w^*$) topology are compact convex
spaces.

 Let $\sigma$ be a stationary measure on the
sequence space $\Lambda^\Z$. The shift operator $S$, defined by
$(Sx)_k=x_{k+1}$, is continuous and the system  $
\mbox{\boldmath$\Lambda$}^\Z= (\Lambda^\Z,\A^\Z,\sigma,S)$ is
measure preserving. We call it the \emph{sequence space system}
determined by the stationary measure $\sigma$. Moreover, we call
the stationary process $\{x_i\}_{i\in\Z}$ the \emph{sequence space
process} determined by $\sigma$. To ease notation we denote by
$x_i$ both the $i$-th coordinate of a point $x$ and the function
that maps each point to its $i$-th coordinate.

Let $\{f_i\}_{i\in\Z}$ be a $\Lambda$-valued stationary process.
If there exists a stationary measure $\sigma$ on $\Lambda^\Z$ that
makes the processes $\{f_i\}_{i\in\Z}$ and $\{x_i\}_{i\in\Z}$
equivalent we say that the second process is the \emph{sequence
space representation} of the first. The next classical result
(\cite{Br}, page 107) is an easy consequence of  Kolmogorov's
extension theorem.

\begin{proposition}\label{P:model1}
Every $\Lambda$-valued stationary process $\{f_i\}_{i\in\Z}$ has a
$\Lambda^\Z$ sequence space representation.
\end{proposition}

Let ${\bf X}$ be a measure preserving system and $\F$ be a
$T$-generating  subalgebra. Suppose that ${\bf X}$ is isomorphic
to a sequence space system ${\bf I}^\Z=(I^\Z,\B^\Z,\sigma,S)$,
where $I=[0,1]$, and the isomorphism $\phi \colon X\to I^\Z$ maps
sets in $\B(\F)$ to $x_0$-measurable sets. Then we say that the
system ${\bf I}^\Z$ is the \emph{sequence space representation of
${\bf X}$ with respect to the subalgebra $\F$}. The next
proposition is a variation of a classical result:

\begin{proposition}\label{P:model2}
Every invertible measure preserving system ${\bf X }$ has a
sequence space representation with respect to any $T$-generating
subalgebra $\F$.
\end{proposition}
\noindent \emph{Sketch of the Proof.} By a classical result of
Rokhlin (\cite{Ro}) the sub-$\sigma$-algebra $\B_0=\B(\F)$ induces
a partition $\P=\{P_t\}_{t\in I}$ of $X$ by $\B_0$-measurable sets
such that every $\B_0$-measurable set is equal (up to a set of
measure zero) to a union of partition elements $P_t$. For every
open set $A\subset I$ we set $A'=\bigcup_{t\in A}P_t$. We define
the measure $\sigma$ on cylinder sets of $I^\Z$ by $
\sigma\big(\bigcap_{i=-n}^n\{x\colon x_i \in A_i \}\big)=
\mu\big(\bigcap_{i=-n}^n\{x\colon T^i x \in A_i'\}\big)$, and then
extended it to the whole sequence space using Kolmogorov's
extension theorem.  We can check that the advertised sequence
space representation of ${\bf X}$ is determined by the  measure
$\sigma$.

\subsection{Van der Corput's Lemma.}\label{L:Van}
The following classical lemma will be needed in the sequel
(\cite{FW}, page 47):
\begin{lemma}[\bf{Van der Corput}]
Let  $\{x_n\}_{n\in\N}$ be a  bounded sequence of vectors in a
Hilbert space. For each $m$ we set
$$
b_m=\overline{\lim}_{N\to\infty}\Big|\frac{1}{N}\sum_{n=1}^N
<x_{n+m},x_n>\Big|.
$$
Assume that
$$
\overline{\lim}_{M\to\infty}\frac{1}{M}\sum_{m=1}^M b_m=0.
$$
 Then
$$
\lim_{N\to\infty}\frac{1}{N}\sum_{n=1}^N x_n=0
$$
in the $L^2$ norm.
\end{lemma}

\section{Strong stationarity}

\subsection{Definitions.} \label{SS:def}
\begin{definition}\label{D:sst} (i) A stationary process
$\{f_i\}_{i\in \Z}$ is \emph{strongly stationary} if the finite
dimensional statistics of $\{f_i\}_{i\in \Z}$ and
$\{f_{ni}\}_{i\in \Z}$ are the same for every $n\in\N$.

(ii) An invertible  measure preserving system ${\bf X}$ is
\emph{strongly stationary} if there exists a $T$-generating
function algebra $\F$ such that every process $\{ \bT^if\}_{i\in
\Z}$ is strongly stationary for $f\in \F$. When we want to also
indicate the subalgebra $\F$ we write that $(\bX,\F)$ is strongly
stationary.
\end{definition}

\begin{remark}
 Equivalently, a system ${\bf X}$ is strongly stationary if there
exists  a $T$-generating function algebra $\F$ such that
$$
\int f_0 \ \bT f_1 \ \cdots \ \bT^k f_k \,d\mu = \int f_0 \ \bT^n
f_1 \ \cdots \ \bT^{kn} f_k \, d\mu
$$
for $n\in \N$, $f_i\in \F$.

\end{remark}
 Let $ \mbox{\boldmath$\Lambda$}^\Z$ be a sequence space system.
The subalgebra $\F_0$ of bounded $x_0$-measurable functions is
$T$-generating. Strong stationarity with respect to $\F_0$ is
equivalent to saying that the maps $\tau_n$ defined by $(\tau_n
x)_i=x_{ni}$ are measure preserving for every $n\in\N$.
\begin{definition}
A measure $\sigma$ on the sequence space $\Lambda^\Z$ is
\emph{strongly stationary} if the system $({\bf \Lambda}^\Z,\F_0)$
is strongly stationary.
\end{definition}

\subsection{Stationary processes majorize strongly stationary
ones.}\label{SS:majorize}


Let $\{f_i\}_{i\in \Z}$ and $\{g_i\}_{i\in \Z}$ be two
$\Lambda$-valued stationary processes with underline measures
$\mu$, $\nu$  correspondingly. We say that $\{f_i\}_{i\in \Z}$
\emph{majorizes} $\{g_i\}_{i\in\Z}$ if
$$
\sup_{n\in\N}\mu\Big(\bigcap_{i=-k}^k\{f_{in}\in A_i\}\Big)\geq
\sup_{n\in\N}\nu\Big(\bigcap_{i=-k}^k\{g_{in}\in A_i\}\Big)
$$
for every $k\in\N$ and open sets $A_i\subset\Lambda$.

Furstenberg and Katznelson (\cite{FK2}) proved that every
stationary process majorizes a strongly stationary one. Actually
they established a much more general result using a strong
selection theorem. The argument given below was suggested by Y.
Peres and gives an easier proof for the case that we are
interested.
\begin{theorem}[{\bf Furstenberg and Katznelson}]\label{T:FK}
Every stationary process majorizes a strongly stationary one.
\end{theorem}
\begin{proof}
By Proposition \ref{P:model1} there exists a stationary measure
$\sigma$ on the sequence space $\Lambda^\Z$ such that the
processes $\{f_i\}_{i\in\Z}$ and $\{x_i\}_{i\in\Z}$ have the same
finite dimensional statistics. For $n\in\N$, let $\tau_n$ be the
map defined on $\Lambda^\Z$ by $(\tau_n x)_i=x_{ni}$. It is
straightforward to check that the measure $\tau_n\sigma$ defined
by $\tau_n\sigma(A)=\sigma(\tau_n^{-1}A)$ is  stationary. If we
denote by $\O$ the closure in the $w^*$-topology of the set of all
convex combinations of the measures $\tau_n \sigma$, $n\in\N$,
then $\O$ is a compact convex subset of $\M_s$. The maps $\tau_n$
commute and act continuously and affinely on $\O$ so by the
Markov-Kakutani fixed point theorem (\cite{Co}, page 151) they
have a common fixed point $\nu\in\O$.

We claim that the stationary process $\{x_i\}_{i\in\Z}$ which is
induced by the measure $\nu$ on the sequence space $\Lambda^\Z$ is
strongly stationary and is majorized by the process
$\{f_i\}_{i\in\Z}$. Indeed, the invariance over each $\tau_n$
proves strong stationarity. Moreover, convex linear combinations
of the measures $\tau_n\sigma$ come arbitrarily close to $\nu$ in
the $w^*$-topology.  Hence, for every $\varepsilon>0$ and choice
of $A_i$'s there exists $n\in\N$ such that
$$
\sigma\big(\bigcap_{i=-k}^k\{x_{in}\in A_i\}\big)\geq
\nu\big(\bigcap_{i=-k}^k\{x_{i}\in A_i\}\big)-\varepsilon.
$$
Since for every $n\in \N$ we have
$\sigma\big(\bigcap_{i=-k}^k\{x_{in}\in A_i\}\big)=$
$\mu\big(\bigcap_{i=-k}^k\{f_{in}\in A_i\}\big)$ and
$\nu\big(\bigcap_{i=-k}^k\{x_i\in A_i\}\big)=$
$\nu\big(\bigcap_{i=-k}^k\{x_{in}\in A_i\}\big)$ the result
follows.
\end{proof}

We  deduce now a similar result for measure preserving systems.

\begin{corollary}
Let ${\bf X}=(X,\B,\mu,T)$ be any (not necessarily invertible)
measure preserving  system and $A$ be a $\B$-measurable set. Then
there exists an invertible strongly stationary system ${\bf
\tilde{X}}= (\tilde{X},\tilde{\B},\tilde{\mu},\tilde{T})$
and $B\in\tilde{\B}$, such that $\mu(A)=\tilde{\mu}(B)$ and
$$
\sup_{n\in\N} \mu\big(\bigcap_{i=0}^k
T^{-in}A\big)\geq\sup_{n\in\N} \tilde{\mu}\big(\bigcap_{i=0}^k
{\tilde{T}}^{-in}B\big)
$$
 for every $k\in\N$.
\end{corollary}
\begin{proof}
Let ${\bf 1}_A$ denote the indicator function of the set $A$.
Using a standard argument we extend the one sided stationary
process $\{f_i\}_{i\in\N}=\{\bT^i {\bf 1}_{A}\}_{i\in\N}$ to a two
sided one. We denote the two sided extension by
$\{f_i\}_{i\in\Z}$. By Theorem \ref{T:FK}  the process
$\{f_i\}_{i\in\Z}$ majorizes a strongly stationary one
$\{g_i\}_{i\in\Z}$. Following the  proof of Theorem \ref{T:FK} we
see that $\{g_i\}_{i\in\Z}=\{\bS^i x_0\}_{i\in\Z}$, where $S$ is
the (invertible) shift transformation on the sequence space
$\{0,1\}^\Z$ with some appropriately  chosen measure $\sigma$. We
let $B=\{x\in\{0,1\}^\Z\colon x_0=1\}$ and set ${\bf \tilde{X}}=
\(\{0,1\}^\Z,\B^\Z,\sigma,S\)$. If  $\F_0$ is  the subalgebra of
bounded $x_0$-measurable functions then $({\bf \tilde{X}},\F_0)$
is strongly stationary. The advertised inequality is valid since
$\{f_i\}_{i\in\Z}$ majorizes $\{g_i\}_{i\in\Z}$. Finally,
following again the proof of Theorem \ref{T:FK} we see that
$\mu(A)=\sigma(B)$.
\end{proof}

\subsection{The $\tau_n$'s and Jenvey's result.} Let ${\bf X}$ be a
strongly stationary system. The next proposition  was proved in
\cite{Je2} and gives useful necessary and sufficient conditions
for strong stationarity:
\begin{proposition}[{\bf Jenvey}] \label{P:t_n}
The measure preserving system  ${\bf X}$ is strongly stationary if
and only if  there exists a $T$-generating function algebra $\F$,
and a family of measure preserving transformations
${\{\tau_n\}}_{n\in \N}$, that leave  every function in $\F$
invariant and such that the operators $\bT$ and
$\mbox{\boldmath$\tau_n$}$ satisfy the commutation relations
\begin{equation}\label{E:t_n}
\mbox{\boldmath$\tau_n$}\bT=\bT^n\mbox{\boldmath$\tau_n$}, \quad n
\in \N.
\end{equation}
 Moreover, we can choose the $\tau_n$'s to satisfy
$\mbox{\boldmath$\tau_{mn}$}=
\mbox{\boldmath$\tau_m$}\mbox{\boldmath$\tau_n$}$, for all $m, n
\in \N$.
\end{proposition}
\begin{remark}
Equivalently, for the point transformations  $T$ and $\tau_n$
relation (\ref{E:t_n}) can be written as
 $(T\tau_n)(x)=(\tau_n T^n)(x)$, for a.e. $x\in X$ and  $n\in \N$.
\end{remark}


Using the multiple weak-mixing theorem (\cite{Fu2}, page 86) it is
easy to see that if a strongly stationary system is weak-mixing
then it is Bernoulli. In \cite{Je2} Jenvey shows that the same
conclusion holds if we just assume ergodicity.
\begin{theorem}[{\bf Jenvey}]\label{T:jenvey}
Every ergodic strongly stationary system is a Bernoulli system.
\end{theorem}

We remark that strong stationarity of a system does not imply that
of its ergodic components, so we cannot use this theorem to
determine the structure of the general strongly stationary system.

Applying the argument in the proof of Theorem~\ref{T:jenvey}  for
the general (not necessarily ergodic) strongly stationary system
we can deduce the following:
\begin{proposition} \label{P:jenvey}
If $\bX$  is a strongly stationary system  and $\lambda$ is  an
eigenvalue of $T$ then $\lambda=1$.
\end{proposition}
Since the proof is too long to reproduce  we just indicate the
strategy. Suppose  that  $\chi$ is a $\lambda$-eigenfunction of
$T$, $\lambda$ is not a root of unity, and that $g\in \V_{-k}^k
T^i \F$. We have
$$
\int \chi \ g\ d\mu=\lim_{N\to\infty}\frac{1}{N}\sum_{n=1}^N \int
\mbox{\boldmath$\tau_n$}\chi\  \mbox{\boldmath$\tau_n$}g \ d\mu.
$$
To show that the limit on the right is zero we use repeatedly Van
der Corput's lemma (Lemma \ref{L:Van}) and relation \eqref{E:t_n}.
As it turns out, it suffices to show that the sequence of
functions
$
u_n=\bT^{-n}(\mbox{\boldmath$\tau_{n+m}$}\chi \
\overline{\mbox{\boldmath$\tau_n$}\chi})
$
converges weakly to zero as $n\to\infty$ for every $m\in \N$. From
relation \eqref{E:t_n} it follows that
$\mbox{\boldmath$\tau_n$}\chi$ is a linear combination of
$\lambda^{1/n}$-eigenfunctions. When $\lambda$ is not a root of
unity  this easily gives that the sequence $u_n$ consists of
"almost" pairwise orthogonal functions and hence converges weakly
to zero. It follows that  $\chi$ is orthogonal to the subalgebra
$\V_{-k}^k T^i \F$ for every $k\in\N$. Since the algebra $\F$ is
$T$-generating we get that $\chi=0$. The case where $\lambda$ is a
nontrivial root of unity is trickier but the strategy of the proof
is similar.

\subsection{Examples of strongly stationary systems.}
\label{S:examples} We give now the basic  examples of strongly
stationary systems.

(i)\emph{ Exchangeable systems.} A  system  is exchangeable if
there exists  a $T$-generating function algebra $\F$ which has
finite dimensional statistics invariant under any permutation of
the time parameter. Bernoulli systems are exchangeable. A theorem
of de Finetti says that a system is exchangeable  if and only if
it is a mixture of Bernoulli systems. We can check that
exchangeability implies strong stationarity but the converse is
not true as the next example shows.

(ii)\emph{ Ergodic components circle rotations.} On the $2$-torus,
with group action addition \!\!$\mod{1}$ and the Haar measure,
define $T(x,y)= (x,y+x)$.
To check that the system  is strongly stationary we use
Proposition \ref{P:t_n}. We let $\F$ be the algebra generated  by
the exponentials in $y$ and define the maps $\tau_n$ by
$\tau_n(x,y)= (nx, y)$. Since $e^{iy}\in\F$ and $e^{ix}=Te^{iy}\
e^{-iy},$ the algebra $\F$ is $T$-generating. Moreover, $\tau_n$
is measure preserving, each $\tau_n$ leaves functions in $\F$
invariant, and
$$
T\tau_n(x,y)=(nx,y+nx)=\tau_n T^n(x,y).
$$
Hence, $\mbox{\boldmath$\tau_n$}\bT=\bT^n
\mbox{\boldmath$\tau_n$}$ for $n\in\N$.

(iii) \emph{ Ergodic components affine transformations on $\T^d$.}
On $\T^3$, with group action addition \!\!$\pmod 1$ and the Haar
measure, define $T(x,y,z)= (x,y+x,z+y)$.  To see that the induced
measure preserving system is strongly stationary we use again
Proposition \ref{P:t_n}. We let $\F$ be the algebra generated by
exponentials in $z$ and define the maps $\tau_n$ by
$\tau_n(x,y,z)= (n^2 x, ny+ \binom{n}{2} x,z)$. We can check as
before that the algebra $\F$ is $T$-generating. Each $\tau_n$
clearly leave functions in $\F$ invariant and a direct computation
shows that it satisfies the right commutation relations.


More generally, on $\T^d$ with group action addition \!\!$\pmod 1$
and the Haar measure,  define
$$
T(x_1,x_2,\ldots,x_d)=(x_1,x_2+x_1,\ldots,x_d+x_{d-1}).
$$
This time $\F$ is the algebra generated by exponentials in $x_d$
and the $\tau_n$'s have the form
$$
\tau_n(x_1,x_2,\ldots,x_d)=\Big(n^{d-1}
x_1,n^{d-2}x_2+\binom{n}{2}x_1,\ldots,nx_{d-1}+\cdots
+\binom{n}{d-1}x_1,x_d\Big).
$$

 (iv)\emph{ Ergodic components affine maps on more general
groups.} In the last  two examples $\T$ can be replaced by any
connected compact abelian group. The form of the $\tau_n$'s and of
the algebra $\F$ is similar. The connectedness assumption is
needed to guarantee that each $\tau_n$ is measure preserving.


Although the previous examples provide an ample supply of strongly
stationary systems,  the building blocks are always Bernoulli
systems and affine transformations on compact abelian groups. In
the last section we will see that this is not the case in general
(Section \ref{SS:last}, example (iv)).

\section{Extremality}

\subsection{Definition of extremality.} \label{SS:extremal}
Recall that by Proposition \ref{P:model2}
every strongly stationary system $({\bf X},\F)$ has a sequence
space representation $(\mathbf{I}^\Z,\F_0)$, where $I=[0,1]$ and
$\F_0$ is the subalgebra of bounded $x_0$-measurable functions.
This representation is completely determined by the measure
$\sigma$ on the sequence space, so characterizing (up to
isomorphism) the strongly stationary measure preserving systems is
equivalent to characterizing the strongly stationary measures on
$I^\Z$. Furthermore, we only have to determine the structure of
the extremal ones, that is, those that cannot be decomposed
nontrivially into a convex combination of strongly stationary
measures. We will make this more precise below.

 Consider the space $\M_{ss}$ of all strongly stationary measures.
 Then  $\M_{ss}$ is a closed convex subset of the space
 of stationary measures $\M_s$ which is $w^{*}$-compact,
 metrizable, and locally convex.
\begin{definition}\label{D:extremal}
(i) Let $\sigma$ be a strongly stationary measure on $I^\Z$. Then
$\sigma$ is \emph{extremal} if it cannot be written in the form
$\sigma=a\sigma_1+(1-a)\sigma_2$ for some $0<a<1$ and strongly
stationary measures $\sigma_1\neq\sigma_2$.

(ii) Let $({\bf X},\F)$ be a strongly stationary system and
$\sigma$ be the strongly stationary measure that determines its
sequence space representation with respect to $\F$. Then $({\bf
X},\F)$ is \emph{extremal} if $\sigma$ is an extremal strongly
stationary measure.
\end{definition}

We will use the  following integral representation theorem of
Choquet (\cite{Ph}):
\begin{theorem}[{\bf Choquet}]  \label{T:choquet}
Let $X$ be a metrizable compact convex subset of a locally convex
space $E$ and $x_0 \in X$. Then there exists a Borel probability
measure $\mu$ on $X$, supported on the extreme points $ext(X)$,
that satisfies $x_0= \int_{ext(X)} x d\mu (x)$ (that is, $l(x_0)=
\int_{ext(X)} l(x) ds(x)$ holds for every $l$ in the dual of $X$).
\end{theorem}

It follows that  the general strongly stationary measure or system
is an integral of extremal ones. So we can focus our attention on
determining the structure of the extremal strongly stationary
measures or systems.


\subsection{Necessary and sufficient conditions.}
It is well known  that the set of extremal points of the space
$\M_s$ is the set of ergodic measures (with respect to the shift
transformation $S$), and that different ergodic measures are
mutually singular. It is not hard to establish the analogous
results for the space of strongly stationary measures $\M_{ss}$.
The corresponding action on $\M_{ss}$ is the joint action of $S$
and the $\tau_n$'s $((\tau_n x)_i=x_{ni})$.
\begin{proposition} \label{P:singular}
Different extremal strongly stationary measures are mutually
singular.
\end{proposition}
\begin{proof}
Let $\mu_1$ and $\mu_2$ be two different extremal strongly
stationary measures and let $S$ denote the shift transformation on
$I^\Z$. Consider the Lebesgue decomposition  of $\mu_1$ with
respect to $\mu_2$, that is, write  $\mu_1= a \nu_1 + (1-a) \nu_2$
where $\nu_1, \nu_2$ are probability measures such that $\nu_1
\perp \mu_2$ and $\nu_2 \ll \mu_2$. Then
\begin{align*}
 \mu_1= S\mu_1 &= a\, S\nu_1+  (1-a)\, S\nu_2, \\
 \mu_1= \tau_n\mu_1 &= a\,\tau_n \nu_1+ (1-a)\,\tau_n \nu_2.
\end
{align*} The Lebesgue decomposition is unique, so  both $S$ and
$\tau_n$  preserve $\nu_1$, $\nu_2 $. This means that $\nu_1$ and
$\nu_2 $ are both strongly stationary measures. Since $\mu_1$ is
extremal we have either $\mu_1 = \nu_1$ or $\mu_1= \nu_2$. If
$\mu_1=\nu_1$ then $\mu_1$ and  $\mu_2$ are mutually singular. So
it remains to show that  $\mu_1\neq\nu_2$. Suppose on the contrary
that $\mu_1=\nu_2$. Then $\mu_1 \ll \mu_2$, so $\mu_1= \int f d
\mu_2$ for some $f \in L^{1}(\mu_2)$. The choice of $f$ is unique,
so we conclude as before that $f$ is $S, \tau_n$-invariant (with
respect to $\mu_2$). Since $f$ is nonconstant ($\mu_1 \neq \mu_2$)
there exists a $S,\tau_n$-invariant set $A$ such that
$0<\mu_2(A)<1$. Then  $\mu_2$ is a nontrivial convex combination
of the induced strongly stationary measures on $A$ and $A^c$. This
contradicts the extremality of $\mu_2$ and completes the proof.
\end{proof}

\begin{proposition} \label{P:extremality} A strongly stationary
measure  $\sigma$ is extremal if and only if the joint action of
$S$ and the $\tau_n$'s is ergodic.
\end{proposition}
\begin{proof}
Suppose that the joint action is not ergodic, that is,  there
exists a $S, \tau_n$-invariant set
 $A$ with $0<\sigma (A)<1$. Then $\sigma$ is a nontrivial convex
 combination of the induced strongly stationary measures on $A$
 and $A^c$. Hence $\sigma$ is not extremal.

Conversely, suppose that the joint action is ergodic. Let
$\sigma=a\sigma_1+(1-a)\sigma_2$, for some strongly stationary
measures  $\sigma_1$, $\sigma_2$. Then $\sigma_1$ is absolutely
continuous with respect to $\sigma$ and  the corresponding
Radon-Nikodym derivative $d\sigma/d\sigma_1$ is a
$S,\tau_n$-invariant function. Since the joint action is ergodic,
$d\sigma/d\sigma_1$ is constant $\mu$-a.e.. Hence
$\sigma=\sigma_1$. This proves that  $\sigma$ is extremal.
\end{proof}

\begin{remarks}
(i) It follows that a strongly stationary system $({\bf X},\F)$ is
extremal if and only if the joint action of $T$ and the $\tau_n$'s
(of Proposition \ref{P:t_n}) is ergodic.

(ii) Using this proposition we can easily check  that the examples
on $\T^d$ given in Section \ref{S:examples} \ \!are extremal.
\end{remarks}

\subsection{Homogeneity property.} In this section we will show that
the ergodic components of an extremal strongly stationary system
enjoy a  homogeneity property, in the sense that their structure
is similar.

\begin{lemma} \label{L:permute}
Let $({\bf X},\F)$ be a strongly stationary system. Then $\tau_n$
leaves the sub-$\sigma$-algebra of $T$-invariant sets invariant.
\end{lemma}
\begin{proof}
Let $A$ be a $T$-invariant set, that is $T^{-1}A=A$. The
commutation relations of Proposition \ref{P:t_n} give that
$$
\tau_n^{-1}A= \tau_n^{-1}T^{-1}A= T^{-n}\tau_n^{-1}A.
$$
Hence, $\tau_n^{-1}A$ is left invariant by $T^n$. By Theorem
\ref{P:jenvey} the transformation  $T$ does not have nontrivial
roots of unity as eigenvalues. It follows that  $\tau_n^{-1}A$ is
also $T$-invariant, completing the proof.
\end{proof}
\begin{lemma} \label{L:full}
Let $({\bf X},\F)$  be an extremal strongly stationary system. If
$A\in \I$ has positive measure then \ \! $\bigcup_{n\in \N}
\tau_n^{-1}(A)$ has full measure.
\end{lemma}
\begin{proof} Let $B= \bigcup_{n\in \N} \,\tau_n^{-1}(A)$.
 In view of  Proposition  \ref{P:extremality} it suffices to show that
$B$ is $T,\tau_n$-invariant. Since $\tau_{mn} = \tau_m \tau_n$ we
have
$$
\tau_m^{-1}(B)=\bigcup_{n\in \N}\tau_m^{-1} \tau_n^{-1}(A)=
\bigcup_{n\in \N}\tau_{mn}^{-1}(A)\subset B.
$$
So the set $B$ is $\tau_m$-invariant for every $m\in \N$.
Moreover, since $T^{-1}B$ is equal to $\bigcup_{n\in \N}
T^{-1}\tau_n^{-1}(A)$ and by Lemma \ref{L:permute} the set
$\tau_n^{-1}(A)$ is $T$-invariant, $B$ is also $T$-invariant.
\end{proof}

\begin{definition}
Let ${\bf X}$  be a measure preserving system  with ergodic
decomposition $\mu=\int \mu_t \ d\lambda (t)$. We say that the
sets $A, B \in \I$ with positive $\lambda$-measure
 are \emph{factor power equivalent} (FPE), if for $\lambda$-a.e.
$b\in B$ there exists
  $ a \in A$ and $n \in \N$ such that  $(X,\B,\mu_b,T)$ is a factor
of $(X,\B,\mu_a,T^n)$, and vice versa.
\end{definition}

We are now ready to prove the advertised homogeneity property.

\begin{theorem} \label{T:key}
Let $({\bf X},\F)$ be  an extremal strongly stationary system with
ergodic decomposition  $\mu=\int \mu_t \ d\lambda(t)$. Then any
two sets  $A,B\in \I$ with positive measure are FPE.
\end{theorem}
\begin{proof}
By Lemma \ref{L:permute} each $\tau_n$ permutes the ergodic
fibers. Suppose that  $\tau_n$ maps the fiber $X_a$ to the fiber
$X_b$ $(b=\tau_n a)$. The pointwise commutation relations $ T
\tau_n = \tau_n T^n$ show that $(X, \B,\mu_b,T)$  is a factor of
$(X, \B,\mu_a,T^n)$. By Lemma \ref{L:full} if $A\in \I$ has
positive measure then $\bigcup_{n\in \N} \tau_n^{-1}(A)$ has full
measure. It follows that  $A$ is FPE to $X$. So any two sets  $A,
B\in \I$ with positive measure are FPE.
\end{proof}

We call a property ``nice'' if it  is preserved by factors and
powers of  measure preserving systems. The homogeneity property
just established  allows us to extend ``nice'' properties from a
nontrivial set of ergodic components to $\lambda$-a.e. ergodic
component.

As an application, suppose that  a nontrivial set (of positive
$\lambda$ measure)  of ergodic components of an extremal strongly
stationary system $(X,\F)$ is weak mixing . We claim that it is a
Bernoulli system. Indeed, weak mixing is a ``nice'' property, so
Theorem \ref{T:key} gives  that $\lambda$-a.e. ergodic component
of the system is weak mixing. Strong stationarity gives
$$
\int f_0 \ \bT f_1 \ \cdots \ \bT^k f_k \
d\mu=\frac{1}{N}\sum_{n=1}^N \int f_0 \ \bT^n f_1 \ \cdots \
\bT^{kn} f_k \ d\mu,
$$
for all $k$, $N\in\N$, $f_i\in \F$. Let $\mu=\int \mu_t \
d\lambda(t)$ be the ergodic decomposition of $\mu$. Letting $N\to
\infty$ and using the multiple weak mixing theorem (\cite{Fu2},
page 86) we get
$$
\int f_0 \ \bT f_1 \ \cdots \ \bT^k f_k \ d\mu=\int \(\int f_0 \
d\mu_t \ \int f_1 \ d\mu_t \ \cdots \ \int f_k \ d\mu_t \)\
d\lambda(t).
$$
It follows that the system is an integral of Bernoulli systems and
being extremal it must be Bernoulli.

\section{Reduction to distal systems} \label{S:results}

\subsection{Characteristic factors.} \label{SS:characteristic}
The notion of a characteristic factor was introduced by
Furstenberg in order to facilitate the study of  nonconventional
ergodic averages. The idea is to find the smallest factor of a
system that completely determines the limit behavior of these
averages and then work with this simpler system.
\begin{definition}
A factor ${\bf Y}=(Y,\E,\nu,T)$ of a system ${\bf X}$ is
\emph{characteristic for} $k$ \emph{terms} if
$$
\lim_{N\to \infty}\frac{1}{N} \sum_{n=1}^N  \Big(\prod_{i=1}^k
\bT^{in}f_i -   \prod_{i=1}^k \bT^{in}\bbE(f_i|\E)\Big)=0
$$
in $L^2({\bf X})$ for  $f_i \in L^{\infty}({\bf X})$.
\end{definition}

Furstenberg   (\cite{Fu1}) proved  for ergodic systems that the
$k$-step distal factor is characteristic for $k+1$ terms. We want
to use  this result for general systems (not necessarily ergodic),
so for completeness we include a proof that covers the general
case.

\begin{theorem}[{\bf Furstenberg}] \label{T:furstenberg}
Let ${\bf X}$ be a measure preserving system. Then the factor
${\bf D}_{k-1}$ is characteristic for $k$-terms.
\end{theorem}
\begin{proof}
We use induction on $k$. For $k=1$ this is the context of the
$L^2$-ergodic theorem. Assume that the statement is valid for $k$,
we will establish it for $k+1$. It suffices to show that if one of
the $f_i$'s is orthogonal to $\D_k$ then
$$
\lim_{N\to \infty} \frac{1}{N} \sum_{n=1}^N \prod_{i=1}^{k+1}
\bT^{in}f_i=0
$$
in $L^2({\bf X})$. Indeed, add and subtract $\bbE(f_i|\D_k)$ to
every $f_i$ and expand the product. All the terms but the two that
we are interested  will converge to zero giving us the desired
identity.

So suppose that $\bbE(f_{1}|\D_k)=0$ (the argument is similar if
$\bbE(f_i|\D_k)=0$ for $i\neq 1$). We apply Van der Corput's lemma
\eqref{L:Van} on the Hilbert space $L^2({\bf X})$ with
$a_n=\prod_{i=1}^{k+1}T^{in}f_i$.
 In order to show that
$$
\lim_{N\to \infty}\frac{1}{N} \sum_{n=1}^N a_n=0
$$
 in $L^2({\bf X})$ it suffices  to establish  that
 \begin{equation}\label{E:A}
\overline{\lim}_{M\to \infty}\frac{1}{M} \sum_{m=1}^M b_m=0,
 \end{equation}
 where
$$
b_m=\overline{\lim}_{N\to\infty} \Big|\frac{1}{N} \sum_{n=1}^N
<a_{n+m},a_n>\Big|= \overline{\lim}_{N\to \infty} \Big|\frac{1}{N}
\sum_{n=1}^N \int \prod_{i=0}^{k} \bT^{in}(\bT^{(i+1)m}f_i
\bar{f}_i)\ d\mu\Big|.
$$
 By the induction hypothesis the last limit  is equal to
$$
\overline{\lim}_{N\to \infty}\Big| \frac{1}{N} \sum_{n=1}^N \int
\prod_{i=0}^{k} \bT^{in}\bbE(\bT^{(i+1)m}f_i \bar{f}_i|\D_{k-1})\
d\mu\Big|.
$$
 Now we make use of the fact that $\bbE(f_{1}|\D_k)=0$. It follows from
Lemma \ref{L:kronecker}  that the function $g(x_1,x_2)=f_1(x_1)
\bar{f_1}(x_2)$ is orthogonal to the space of invariant functions
of ${\bf X\times_{D_{k-1}} \! X}$.
If $S=T\times T$,  applying the $L^2$-ergodic theorem for the
system  ${\bf X\times_{D_{k-1}} \! X}$ we get
\begin{equation}\label{786}
 \lim_{M\to \infty}\frac{1}{M} \sum_{m=1}^M \int \bS^mg
\  \bar{g} \ d \mu'=
 \lim_{M\to \infty}\frac{1}{M} \sum_{m=1}^M \int \left| \bbE(\bT^m
f_1\bar{f}_1|\D_{k-1}) \right|^2  d\mu=0.
\end{equation}
Since  every $f_i$ is bounded the Cauchy-Schwartz inequality gives
$$
\Big|\int \prod_{i=0}^{k} \bT^{in}\bbE(\bT^{(i+1)m}f_i
\bar{f}_i|\D_{k-1})\ d\mu\Big|\leq  L\int \left| \bbE(\bT^m
f_1\bar{f}_1|\D_{k-1}) \right|^2  d\mu
$$
for some number $L$ that is independent of $n$. Hence,
$$
b_m \leq L\int \left| \bbE(\bT^m f_1\bar{f}_1|\D_{k-1}) \right|^2
d\mu.
$$
From this and (\ref{786})  it follows that the limit in
(\ref{E:A}) is $0$. This completes the induction.
\end{proof}

\begin{corollary}\label{C:D}
The   distal factor of a system is a characteristic factor for
$k$-terms for every $k\in\N$.
\end{corollary}


\subsection{Relative Bernoulli extensions.}
The notion of a relative Bernoulli extension was introduced in
\cite{Th2}.
\begin{definition}
Let $\bX$  be an ergodic system and $\bY$ be a factor of $\bX$.
Then $\bX$  is a \emph{relative Bernoulli} extension of $\bY$ if
$\bX$ is isomorphic to the direct product of a Bernoulli system
and ${\bf Y}$.
\end{definition}

Let $(X,\B,\mu)$ be a measure space. Two \emph{finite}
$\B$-measurable partitions $\P=\{P_i\}_{i=1}^k$ and
$\mathcal{Q}=\{Q_i\}_{i=1}^l$ of $X$ are
\emph{$\varepsilon$-independent} if
$$
\sum_{i,j} \left|\mu(P_i\cap Q_j)-\mu(P_i)\mu(Q_j)\right|\leq
\varepsilon.
$$

Let $\bX$ be an invertible and ergodic system, $\bY$ be a factor
of $\bX$, and $\mu=\int \mu_y d\nu(y)$  be the disintegration of
$\mu$ over $\bY$. A sequence of finite partitions $\{\P_i\}_{i\in
\Z}$ is \emph{weak Bernoulli relative to} $\bY$ if for $\nu$-a.e.
$y$ the following is true:
for given $\varepsilon>0$ there exists
$N\in\N$ such that for every $m\geq 1$
$$
\V_{-m}^0 \P_i \quad {\rm is \ } \e {\rm-independent \ of } \quad
\V_{N}^{N+m}\P_i.
$$
with respect to $\mu_y$.

The following theorem is a consequence of the results of the
articles \cite{Th1} (Propositions 3, 4, and 5) and \cite{Th2}
(Lemma 6). One can also deduce this  from Theorem $2$ in
\cite{Ki}.

\begin{theorem}[{\bf Thouvenot}]\label{T:thouvenot}
Let $\bX$ be an invertible, ergodic system, and  $\bY$ be a factor
of $\bX$. Suppose  that for some finite $T$-generating partition
$\P$ the sequence of partitions $\{T^i\P\}_{i\in\Z}$
 is weak Bernoulli relative to $\bY$. Then $\bX$ is a relative Bernoulli
extension of $\bY$.
\end{theorem}
 Note that the relative notion of weak Bernoulli is a stronger property
 than the relative notion of very weak Bernoulli  that was used in \cite{Th2}.

\subsection{The relative Bernoulli property.}
 \label{SS:rwm}
\begin{lemma}\label{L:distal}
Let $({\bf X},\F)$ be a strongly stationary system and
$\{\tau_n\}_{n\in\N}$ be the maps defined in Proposition
\ref{P:t_n}. Then the spaces $L^2(\bD_k)$ and $L^2({\bf D})$ are
$\tau_n$-invariant for every $n\in \N$.
\end{lemma}
\begin{proof}
Let $n\!\in \!\N$. By Lemma \ref{L:permute} the subspace $L^2({\bf
D}_0)$ is $\tau_n$-invariant. From the definition of ${\bf D}_k$
and ${\bf D}$ it suffices  to show that every maximal compact
extension of a $\tau_n$-invariant space is also
$\tau_n$-invariant. So suppose that ${\bf X}$ is a maximal compact
extension of ${\bf Y}$ and that $L^2({\bf Y})$ is
$\tau_n$-invariant. Consider the disintegration $\mu= \int \mu_y
d\nu(y)$ of $\mu$ over ${\bf Y}$. It suffices to show that if  $f$
is compact relative to ${\bf Y}$ then so is $\tau_n f$. Let $\e>
0$. There exists a finite set of functions $g_1,\ldots,g_m$ such
that
$$
\min_{1\leq s\leq m} \|\bT^i f- g_s\|_{L^2(\mu_y)} <\e
$$
for every $i\in \N$, and $\nu$-a.e. $y\in Y$. Write $i=i' n + r$,
for some $i'\in \N$ and $0\leq r\leq n-1$. Using the commutation
relations of Proposition \ref{P:t_n} we get
$$
f(\tau_nT^ix)=f(T^{i'}\tau_nT^rx).
$$
So
$$
\|\bT^i(\mbox{\boldmath$\tau_n$} f)-\bT^r(\mbox{\boldmath$\tau_n$}
 g)\|_{L^2(\mu_y)}=\| f(T^{i'}\tau_nT^rx)-
g(\tau_nT^rx)\|_{L^2(\mu_{y})}= \|\bT^{i'}f-g\|_{L^2(\mu_{y'})},
$$
where $y'=\tau_n T^r y$. The last equality is valid since
$L^2({\bf Y})$ is  invariant under both $\tau_n$ and $T$. It
follows that  the set of functions
$\left\{\bT^r(\mbox{\boldmath$\tau_n$}  g_s), \  1\leq s\leq m,\
0\leq r\leq n-1\right\}$ is fiberwise a finite  $\e$-net relative
to ${\bf Y}$ for the orbit $\left\{\bT^i(
\mbox{\boldmath$\tau_n$}f) \right\}_{i\in \N}$. This shows that
the function $ \mbox{\boldmath$\tau_n$}f$ is compact relative to
{\bf Y} and completes the proof.
\end{proof}
\begin{remark}
Since we are only going to use  the $\tau_n$-invariance of
$L^2(\bD)$ we can avoid the use of Lemma \ref{L:permute}. Indeed,
the distal factor can be exhausted by a sequence (possibly
transfinite) of maximal isometric extensions starting from the
trivial factor (determined by the algebra of constant functions).
We can then use the step by step argument of the previous proof to
show that $L^2(\bD)$ is $\tau_n$-invariant.
\end{remark}

\begin{theorem}\label{T:rwm}
Let $({\bf X},\F)$  be a strongly stationary system. Then its
distal factor  is strongly stationary and almost every ergodic
component of ${\bf X}$ is a relative Bernoulli extension of its
distal factor.
\end{theorem}
\begin{proof}
{\bf Step 1.} Let ${\bf D}=(D,\D,\nu,T)$ be the distal factor of
${\bf X}$.
 Strong stationarity gives
$$
\int{ f_0 \, \bT f_1 \, \cdots \ \bT^{k} f_k}\, d\mu = \int{ f_0
\, \bT^n f_1 \, \cdots \  \bT^{kn} f_k}\, d\mu
$$
for all $k,n\in\N$, $f_i\in \F$. Averaging over $n$ and taking the
limit as $N\to$ $\infty$ gives
$$
 \int{ f_0 \, \bT f_1 \, \cdots \
\bT^{k} f_k}\, d\mu = \lim_{N\to\infty}\frac{1}{N} \sum_{n=1}^N
\int{ f_0 \, \bT^n f_1 \, \cdots \  \bT^{kn} f_k} \,d\mu.
$$
By Theorem \ref{T:furstenberg} the last average is equal to
\begin{equation}\label{E:average}
\lim_{N\to\infty}\frac{1}{N} \sum_{n=1}^N \int{ \bbE(f_0|\D) \,
\bT^n \bbE(f_1|\D) \, \cdots \ \bT^{kn} \bbE(f_k|\D)} \,d\nu.
\end{equation}
The maps $\{\tau_n\}_{n\in\N}$ leave the functions in $\F$
invariant, as well as the space  $L^2({\bf D})$ (by Lemma
\ref{L:distal}), hence
\begin{equation}\label{E:t_n2}
\mbox{\boldmath$\tau_n$}
\bbE(f_i|\D)=\bbE(\mbox{\boldmath$\tau_n$} f_i|\D)=\bbE(f_i|\D).
\end{equation}
Since $\mbox{\boldmath$\tau_n$}\bT=\bT^n\mbox{\boldmath$\tau_n$} $
and $\tau_n$ is measure preserving for $n\in\N$,  we get using
\eqref{E:t_n2} that all the integrals in (\ref{E:average}) are
equal to
$$
\int{ \bbE(f_0|\D) \, \bT \bbE(f_1|\D) \, \cdots \ \bT^k
\bbE(f_k|\D)} \,d\nu.
$$
 This shows that
\begin{equation} \label{E:distal}
\int{ f_0 \, \bT f_1 \, \cdots \ \bT^{k} f_k} \,d\mu = \int{
\bbE(f_0|\D) \, \bT \bbE(f_1|\D) \, \cdots \ \bT^k \bbE(f_k|\D)}
\,d\nu
\end{equation}
for  $k\in \N$, $f_i\in \F$.

{\bf Step 2.} We will strengthen (\ref{E:distal}) to a  fiberwise
relation and prove the first claim. Call $\F_\D$ the algebra
generated by functions of the form $\bbE(f|\D)$, where $f\in \F$.
Let $\D'$ be the sub-$\sigma$-algebra of $\V_{i=0}^\infty
\bT^{i}\F_\D$-measurable sets. Clearly $\D'\subset \D$. Since
$\D'$ is $T$-invariant it induces a factor ${\bf
D}'=(D',\D',\nu',T)$ of ${\bf X}$. Applying $\tau_n$ to the left
hand side of the equation below and using the previous averaging
technique we get as before that
\begin{equation} \label{E:Long}
\int \prod_{i=0}^k\bT^if_i  \ \prod_{j=0}^m \bT^j g_j\ d\mu= \int
\prod_{i=0}^k \bT^i\bbE(f_i|\D)\ \prod_{j=0}^m \bT^jg_j\ d\nu
\end{equation}
for  $k,m\in \N$, $f_i\in \F$, $g_i \in\F_\D$. Observe that
$\bbE(f|\D)$  is $\D'$-measurable for $f \in \F$, so $\D$ can be
replaced by  $\D'$ in (\ref{E:Long}).  Moreover, all the $g_i$'s
are $\D'$-measurable, so (\ref{E:Long})  takes the following form
\begin{equation}
\int \Big[ \bbE\Big(\prod_{i=0}^k\bT^if_i|\D'\Big)-  \prod_{i=0}^k
\bT^i\bbE(f_i|\D')\Big] \prod_{j=0}^m \bT^j g_j\ d\nu'= 0.
\end{equation}
Since $\V_{i=0}^\infty T^i \F_\D$ is dense in $L^2({\bf D}')$ we
get
$$
\int  \Big[ \bbE\Big(\prod_{i=0}^k\bT^if_i|\D'\Big)-
\prod_{i=0}^k \bT^i\bbE(f_i|\D')\Big] \ g \ d\nu'=0
$$
for every $g\in L^2({\bf D'})$. This can only happen if
\begin{equation}\label{E:D'}
 \bbE\Big(\prod_{i=0}^k\bT^if_i|\D'\Big)(y)=\prod_{i=0}^k \bT^i\bbE(f_i|\D')(y)
\end{equation}
for $\nu'$-a.e. $y\in D'$.

Next  we  claim that $\D'=\D$. Relation  (\ref{E:D'}) easily
implies that  ${\bf X}$ is a  relative weak mixing extension of
${\bf D}'$. Since {\bf D} is the minimal factor with respect to
which ${\bf X}$ is relative weak mixing, $\D$ must be contained in
$\D'$. Thus, $\D'=\D$ and (\ref{E:D'}) takes the form
\begin{equation}\label{E:D}
 \bbE\Big(\prod_{i=0}^k\bT^if_i|\D\Big)(y)=\prod_{i=0}^k \bT^i\bbE(f_i|\D)(y)
\end{equation}
for $\nu$-a.e. $y\in D$.

Since every $f\in \F_\D$ is $\tau_n$-invariant and $\F_\D$ is a
$T$-generating algebra for $\D$, the system {\bf D} is strongly
stationary.

{\bf Step 3.} We will now prove the second claim. First assume
that the sub-$\sigma$-algebra $\B(\F)$ is determined by a finite
partition $\P$. Then $\P$ is $T$-generating for almost every
ergodic component and relation (\ref{E:D}) is valid for almost
every ergodic component, provided that we replace $\D$ with the
distal factor of the corresponding ergodic component. To simplify
the notation we assume that ${\bf X }$ is ergodic, and we keep in
mind that any result  we get will be valid for almost every
ergodic component.

We claim that the sequence of partitions $\{T^i\P\}_{i\in\Z}$
satisfies the conditions of  Theorem \ref{T:thouvenot}. To see
this  observe  first that if we replace $k$ with $2k$ in
(\ref{E:D}) and then apply $\bT^{-k}$  we get
\begin{equation}\label{E:D-}
 \bbE\Big(\prod_{i=-k}^k\bT^if_i|\D\Big)(y)=\prod_{i=-k}^k
 \bT^i\bbE(f_i|\D)(y)
\end{equation}
for every $k\in \N$, $\nu$-a.e. $y\in Y$, and $\P$-measurable
functions $f_i$. Hence,
$$
\int f g \ d\mu_y=\int f\ d\mu_y \ \int g \ d\mu_y
$$
for every $k\in\N$, and $\nu$-a.e. y,  whenever $f$ is
$\V_{-k}^0T^i\P$-measurable and $g$ is $\V_{1}^kT^i\P$-measurable.
So $\{T^i\P\}_{i\in \Z}$ is weak Bernoulli with respect to $\mu_y$
for $\nu$-a.e. $y$. Theorem \ref{T:thouvenot} now implies that the
system is a relative Bernoulli extension of its distal factor.
Hence, almost every ergodic component is isomorphic to the
 direct product of its distal factor and a Bernoulli system.

In general, since $L^2({\bf X})$ is separable there exists a
sequence $\{\P_n\}_{n\in\N}$ of finite $\F$-measurable partitions
such that $\P_{n+1}$ refines $\P_n$ and
$\V_{n=1}^\infty\P_n=\B(\F)$. Applying the previous argument for
the factor that is $T$-generated by $\P_n$ we get a factor ${\bf
X}_n$ which is isomorphic to the direct product of a distal system
and a Bernoulli system. Since increasing unions of distal systems
is distal and of Bernoulli systems  Bernoulli (\cite{Or}, page
$52$), it follows that almost every ergodic component of ${\bf X}$
is isomorphic to the direct product of a distal system and a
Bernoulli system. This completes the proof.
\end{proof}
\subsection{The distal factor coincides with the $\C$-factor.}
\label{SS:cfactor}

 Let {\bf X} be a measure preserving system. For fixed $k\in \N$
consider the closed subalgebra generated by weak limits (in
$L^2({\bf X})$) of the form
\begin{equation}\label{E:C_k}
\text{w}\!\!\lim_{l\to\infty}\frac{1}{N_l} \sum_{n=1}^{N_l}
 \, \bT^{a_1 n} f_1 \,\,  \bT^{a_2 n} f_2 \, \cdots \, \bT^{a_k n} f_k,
\end{equation}
where we are free to  choose any  $a_1,\ldots,a_k \in \N$,
functions $f_1,\ldots,f_k\in L^{\infty}({\bf X})$, and  an
increasing sequence $\{N_l\}_{l\in \N}$ of positive integers that
guarantees weak convergence. This algebra is conjugation closed
and $T$-invariant so it gives rise to a factor  ${\bf C}_k$. If in
addition we are free to choose any $k\in \N$ we get a factor that
extends every  ${\bf C}_k$. We denote it by ${\bf C}$ and call it
the $\C$-\emph{factor} or \emph{characteristic factor} of ${\bf
X}$. Using the notation of Section \ref{SS:characteristic} the
factor ${\bf C}_k$ is characteristic for $k$ terms and the factor
${\bf C}$ is characteristic for any number of terms. Moreover, by
Theorem \ref{T:furstenberg} we have ${\bf C}_k\subset {\bf
D}_{k-1}$ for every $k\in\N$ and so ${\bf C}\subset{\bf D}$.


First, we prove that the $\C$-factor is $\tau_n$-invariant for
every $n\in\N$.
\begin{lemma} \label{L:invariance}
Let  $({\bf X},\F)$ be a strongly stationary system and
$\{\tau_m\}_{m\in\N}$ be the maps defined in Proposition
\ref{P:t_n}. Then for all $k,m\in\N$ the subspaces  $L^2({\bf
C}_k)$ and $L^2({\bf C})$ are $\tau_m$-invariant.
\end{lemma}
\begin{proof}
Let $k, m\in \N$. Using the commutation relations
$\mbox{\boldmath$\tau_m$}\bT^n= \bT^{nm}\mbox{\boldmath$\tau_m$}$
(see Proposition \ref{P:t_n}) we see that  the operator
$\mbox{\boldmath$\tau_m$}$ maps functions of the form
(\ref{E:C_k}) to functions of the form
$$
\text{w}\!\!\lim_{l\to\infty}\frac{1}{N_l} \sum_{n=1}^{N_l} \
\bT^{a_1n m}(\mbox{\boldmath$\tau_m$} f_1)\ \  \bT^{a_2n
m}(\mbox{\boldmath$\tau_m$}f_2)\, \cdots \ \bT^{a_kn m}
(\mbox{\boldmath$\tau_m$}f_k),
$$
belongs  again to $L^2({\bf C}_k)$. Since functions  of the form
(\ref{E:C_k}) generate $L^2({\bf C}_k)$ the subspace $L^2({\bf
C}_k)$ is $\tau_m$-invariant. A similar argument applies for
$L^2({\bf C})$.
\end{proof}

\begin{theorem}\label{T:C}
Every zero entropy (and hence every distal) strongly stationary
system coincides with its $\C$-factor.
\end{theorem}
\begin{proof}
 Since $\C$ is $T$-invariant, and $\F$ is $T$-generating, it
 suffices to show that every $f\in \F$ is \ $\C$-measurable.
Equivalently,  if $f\in\F$ and $g=f- \bbE(f|\C)$,   we need to
show that $g=0 $. By Lemma \ref{L:invariance} the subspace
$L^2({\bf C})$ is $\tau_n$-invariant. Moreover, $f$ is
$\tau_n$-invariant so $g$ is $\tau_n$-invariant for $n\in\N$.
Given functions $f_1,\ldots,f_k\in \F$, there exists an increasing
sequence $\{N_l\}_{l\in\N}$ of positive integers such that the
sequence
$$
\frac{1}{N_l} \sum_{n=1}^{N_l} \bT^n f_1 \, \cdots \, \bT^{kn} f_k
$$
converges weakly in $L^2({\bf X})$ as $m\to\infty$. Since
$\mbox{\boldmath$\tau_n$}\bT=\bT^n\mbox{\boldmath$\tau_n$}$, the
functions $f_i$ and $g$ are $\tau_n$-invariant,  and $\tau_n$ is
measure preserving for $n\in\N$ we get
$$
\int{ g \, \bT f_1 \, \cdots\, \bT^k f_k} \,d\mu=
\lim_{l\to\infty}\frac{1}{N_l} \sum_{n=1}^{N_l} \int g\, \bT^n f_1
\, \cdots \, \bT^{kn} f_k \, d\mu
$$
for  $k\in\N$,  $f_i\in \F$. Since $\bbE(g |\C) = 0$, we have
$\bbE(g|\C_k) = 0$ and by the definition of $\C_k$ the last
average is $0$. Thus, $g$ is orthogonal to the closed subalgebra
spanned by bounded $\V_{1}^\infty T^{i}\F$-measurable functions.
Since ${\bf X}$  has zero entropy we have  $\V_{1}^\infty
T^{i}\F=\V_{-\infty}^\infty T^{i}\F$. It follows that  $g$ is
orthogonal to the full algebra of the system. Hence $g = 0$,
proving our claim.
\end{proof}



\section{Distal strongly stationary systems}

\subsection{Nilrotations.}\label{SS:nilsystems} Let $G$ be a locally compact
and  separable Lie group. The commutator of  two elements $g,h\in
G$ is the element $[g,h]=g^{-1}h^{-1}gh$. If $A,B\subset G$ we
write $[A,B]$ for the subgroup generated by $\{[a,b]\colon a\in A,
b\in B\}$. The lower central series of $G$ is defined as follows,
$G^{(0)}=G$, $G^{(i)}=[G,G^{(i-1)}]$.  The group $G$ is
\emph{nilpotent of order $k$} if $G^{(k)}=\{e\}$, where $e$
denotes the identity element of $G$. With $G_0$ we denote the
connected component of the identity element of $G$. If $\gG$ is a
discrete subgroup (not necessarily normal) of an order $k$
nilpotent group such that $G/\gG$ is compact we call $G/\gG$ an
\emph{order $k$ nilmanifold}. The group $G$ acts on $G/\gG$ by
left translation $T_a(x\gG)=(ax)\gG$. There exists a unique
probability measure on $G/\gG$ that is invariant under left
translations, we denote it by $\mu$ and call it the \emph{Haar
measure} on ${\bf G}/\Gamma$. If $G$ is nilpotent of order $k$, we
call the system ${\bf G}/\gG=(G/\gG,\mathcal{G}/\gG,T_a,\mu)$ an
\emph{order $k$ nilsystem} and the transformation $T_a$ an
 \emph{order $k$ nilrotation}. An inverse limit of (order $k$)
 nilsystems (nilmanifolds) is called an (order $k$) \emph{pro-nilsystem}
 (\emph{pro-nilmanifold}).

The following generalization of a theorem of Parry (\cite{Pa}) was
proved by Leibman \cite{Lei2}. It was originally established  by
Lesigne \cite{Le3} under an extra hypothesis. Bergelson and Host
(\cite{BH}) showed that this extra hypothesis can be removed, thus
providing another independent proof.

\begin{theorem}[{Leibman}] \label{T:ergodic}
Let $(G/\gG,T_a)$ be a nilsystem and set $Z=G/G^{(1)}\Gamma$. If
$G$ is spanned by the connected component $G_0$ and $a$ then


 (i) The nilsystem $(G/\gG,T_a)$ is uniquely ergodic if and only if the
 factor system  $(Z,T_a)$ is ergodic.

 (ii) If $(G/\gG, T_a)$ is ergodic then its Kronecker factor is
 $(Z,T_a)$.
\end{theorem}
\begin{remark}
As it was noted in \cite{BH} if the nilsystem $(G/\gG,T_a)$ is
ergodic then the projection of $<G_0,a>$ on $G/\gG$ being an open
invariant set is equal to $G/\gG$. Hence,
$G/\gG=<G_0,a>/(\gG\cap<G_0,a>)$. So if $T_a$ is ergodic we can
assume that $G=<G_0,a>$.
\end{remark}

\noindent {\bf Examples.} (i) The prototypical example of a
nonabelian order two ergodic nilsystem is defined on the
\emph{Heisenberg} nilmanifold. Let
$$ \Small{G= \left\{
\begin{pmatrix}
1&x_1&x_3\\0&1&x_2\\0&0&1
\end{pmatrix},  \
x_i\in \R \right\}, \
 \Gamma= \left\{
\begin{pmatrix}
1&k_1&k_3\\0&1&k_2\\0&0&1
\end{pmatrix}, \
 k_i\in\Z \right\} , \ a= \begin{pmatrix}
1&a_1&a_3\\0&1&a_2\\0&0&1
\end{pmatrix}},
$$ where $a_1,a_2\in\R$ are rationally independent and $a_3\in\R$.
 Then $G$ with the standard metric is locally compact and connected.
 Moreover, if the group action is matrix multiplication  then $G$
 is nilpotent of order two
 and $G/\gG$ is compact. So $T_a$ defines a nilsystem  on $ G/\gG$.
 Observe that  $G/G^{(1)
 }\gG\simeq \T^2$ and that the rotation
 on $\T^2$ by $(a_1,a_2)$ is ergodic.  It follows from Theorem \ref{T:ergodic}
 that $(G/\gG, T_a)$ is uniquely ergodic and its Kronecker
 factor is induced by the functions on $x_1,x_2$.

(ii) Let
$$ \Small{G= \left\{
\begin{pmatrix}
1&k&x_3\\0&1&x_2\\0&0&1
\end{pmatrix}\!,
\Small{k \! \in \! \Z, \ \! x_i \! \in \! \R} \right\}, \
 \Gamma= \! \left\{
\begin{pmatrix}
1&k_1&k_3\\0&1&k_2\\0&0&1
\end{pmatrix}\!,
 k_i\in\Z \right\} , \ a= \begin{pmatrix}
1&1&0\\0&1&b\\0&0&1
\end{pmatrix}},
$$ where $b$ is irrational. Then $(G/\gG,T_a)$ is uniquely ergodic and  isomorphic to
the affine system on $\T^2$ with the Haar measure induced by
$T(x,y)=(x+b,y+x)$.

It turns out that every  measure preserving system that is induced
by some distal affine transformation on  compact abelian group
with the Haar measure is isomorphic to  a nilsystem.
But not every nilsystem is isomorphic to an affine system. For
example the order two nilsystem of example (i)  is not
(\cite{Fu3}, page 52).

\subsection{Nonconventional ergodic averages.}
The following  theorem of Host and Kra (\cite{HK3}) is crucial for
our study.
\begin{theorem}[\bf{Host and Kra}]\label{T:Bryna}
Let ${\bf X}$ be an invertible ergodic measure preserving system.
Then the averages
$$
\lim_{N\to\infty}  \frac{1}{N} \sum_{n=1}^N  \ \bT^n f_1\ \bT^{2n}
f_2 \ \cdots \ \bT^{kn}f_k
$$
converge in $L^{2}({\bf X})$ for $f_i\in L^{\infty}({\bf X})$.
Moreover, the characteristic factor (defined in Section
\ref{SS:cfactor}) for these averages is an order $k-1$
pro-nilsystem.
\end{theorem}

If $T$ is an ergodic order $k$ nilrotation it is possible to find
an explicit formula for the limit. This was done for $G$ connected
and $k=3$ by Lesigne (\cite{Le1}) and for general (not necessarily
connected) $G$ and $k$ by Ziegler (\cite{Zi}). To describe the
formula it is convenient to first establish some notation.

Let ${\bf G}/\gG$ be an order $l$ ergodic nilsystem.
It turns out (\cite{Lei1}) that for every $k\in \N$ the set
$$
H_k=\Big\{\big(x_1,\ x_1^2x_2,\ \ldots,\
x_1^{\binom{k}{1}}x_2^{\binom{k}{2}} \cdots \
x_l^{\binom{k}{l}}\big),\  x_i\in G^{(i-1)}\Big\}
$$
is a closed subgroup of $G\times\cdots\times G$ (the product has
$k$ terms) with group action coordinatewise  multiplication. If
$\Delta_k=H_k\cap(\gG\times\cdots\times\gG)$, then the quotient
$H_k/\Delta_k$ is again a nilmanifold and supports a unique left
invariant (under translations by elements in $H_k$) measure that
we denote by $\nu_{H_k}$.
\begin{theorem}[\bf{Ziegler}] \label{T:Ziegler}
Let $(G/\gG,T_a)$ be an order $l$ ergodic nilsystem  (assume that
 $G=<G_0,a>$). If $f_1,\ldots,f_k\in L^\infty( G/\gG)$ then for
almost every $x\in G$ we have
\begin{align*}
\lim_{N\to\infty}  \frac{1}{N} \sum_{n=1}^N &
 \ \bT_a^n f_1(x\gG) \
\cdots  \  \bT_a^{kn}f_k(x\gG)=\\
& \int_{H_k/\Delta_k} \! f_1(xy_1\gG)\ \cdots \ f_k(xy_k\gG) \
d\nu_{H_k}(y\Delta_k)
\end{align*}
where $y=(y_1,\ldots,y_k)$.
\end{theorem}

The next identity will enable us later to  give a general method
for constructing strongly stationary systems starting from totally
ergodic measure preserving systems. It is a consequence of
Theorems \ref{T:Bryna} and \ref{T:Ziegler}.

\begin{theorem} \label{C:Ziegler}
Let ${\bf X}$ be an ergodic measure preserving system such that
$T^r$ is ergodic for some $r\in\N$. Then
\begin{equation}\label{E:Ziegler}
\lim_{N\to\infty}  \frac{1}{N} \sum_{n=1}^N \int f_0 \ \bT^n f_1
\cdots \bT^{kn} f_k \ d\mu=\lim_{N\to\infty}  \frac{1}{N}
\sum_{n=1}^N \int f_0 \ \bT^{rn} f_1 \cdots \bT^{krn} f_k \ d\mu
\end{equation}
for  $f_i\in L^{\infty}({\bf X})$.
\end{theorem}
\begin{proof}
Assume first that $T_a$ is an order $l$ nilrotation defined on
$G/\gG$. Since $T_a^r$ is ergodic $T_a$ is also ergodic, so
Theorem \ref{T:Ziegler} gives that
\begin{align*}
\lim_{N\to\infty}  \frac{1}{N} \sum_{n=1}^N &\int_{G/\gG}
f_0(x\gG) \ \bT_a^n
f_1(x\gG) \ \cdots  \  \bT_a^{kn}f_k(x\gG)\ d \mu(x\gG)=\\
& \int_{G/\gG}\int_{H_k/\Delta_k} f_0(x\gG)f_1(xy_1\gG)\ \cdots \
f_k(xy_k\gG) \ d\nu_{H_k}(y \Delta_k) \  d\mu(x\gG),
\end{align*}
where $f_i\in L^{\infty}({\bf N}/\gG)$, and $H_k$, $\Delta_k$,
$\nu_{H_k}$ are defined as before.
Observe that
the integral on the right does not depend on the nilrotation $T_a$
on $G/\gG$ as long as it is ergodic. Since $T_{ra}=T_a^r$ is
assumed to be ergodic we get that
\begin{align*}
\lim_{N\to\infty} & \frac{1}{N} \sum_{n=1}^N \int_{G/\gG} f_0 \
\bT_a^n
f_1(x\gG)\  \cdots  \ \bT_a^{kn}f_k(x\gG) \ d\mu(x\gG) = \\
\lim_{N\to\infty} &\frac{1}{N} \sum_{n=1}^N \int_{G/\gG} f_0 \
\bT_{ra}^n
f_1(x\gG)\  \cdots \  \bT_{ra}^{kn} f_k(x\gG) \ d\mu(x\gG) =\\
\lim_{N\to\infty}&  \frac{1}{N} \sum_{n=1}^N \int_{G/\gG} f_0 \
\bT_a^{rn} f_1(x\gG)\  \cdots \ \bT_a^{rkn}f_k(x\gG)\ d\mu(x\gG).
\end{align*}
So (\ref{E:Ziegler}) holds when $T$ is a nilrotation. A standard
approximation argument  shows that (\ref{E:Ziegler}) holds when
$T$ defines a pro-nilsystem.  In the general case we have
\begin{align*}
\lim_{N\to\infty}&  \frac{1}{N} \sum_{n=1}^N \int f_0 \ \bT^n f_1
\ \cdots \ \bT^{kn} f_k \ d\mu=\\
 \lim_{N\to\infty}&  \frac{1}{N}
\sum_{n=1}^N \int \bbE(f_0|\C) \ \bT^{n}\bbE(f_1|\C)\  \cdots \
\bT^{kn} \bbE(f_k|\C) \ d\mu,
\end{align*}
where $\C$ is the characteristic factor of the system. It follows
from Theorem \ref{T:Bryna} that the factor $\C$ is a
pro-nilsystem, so the limit on the right is equal to
$$
 \lim_{N\to\infty} \frac{1}{N}
\sum_{n=1}^N \int \bbE(f_0|\C) \ \bT^{rn}\bbE(f_1|\C)\  \cdots \
\bT^{rkn} \bbE(f_k|\C) \ d\mu.
$$
Since the factor $\C$ is characteristic for all terms  the last
limit is equal to
$$
 \lim_{N\to\infty}  \frac{1}{N}
\sum_{n=1}^N \int f_0 \ \bT^{rn}f_1\  \cdots \ \bT^{rkn} f_k \
d\mu.
$$
The result follows.
\end{proof}
We remark that  identity (\ref{E:Ziegler}) was proved for $k=3$ by
Host and Kra (\cite{HK2}). It was also shown there how it implies
an odd version of the ergodic Szemer\'edi theorem. For general $k$
we get:

\begin{corollary}
Let ${\bf X}$ be a measure preserving system. Assume that  $T^r$
is ergodic for some $r\in\N$ and let $A$ be a measurable set with
$\mu(A)>0$. Then for every $0\leq j<r$ we have
$$
\mu\big( A\cap T^{-n}A\cap\cdots\cap T^{-kn}A\big)>0
$$
for some $n\equiv j\pmod{r}$.
\end{corollary}
The odd Szemer\'edi theorem corresponds to the case $r=2$, $j=1$.

\subsection{The structure of the ergodic components.}
Theorem \ref{T:Bryna} combined with Theorems \ref{T:rwm} and
\ref{T:C} enables us to determine the structure of the ergodic
components of the general strongly stationary system.

\begin{theorem} \label{T:components}
Almost every  ergodic component of a strongly stationary system is
isomorphic to the direct product of a Bernoulli system and a
totally ergodic pro-nilsystem.
\end{theorem}
\begin{proof}
Theorem \ref{T:rwm} shows that almost every ergodic component of a
strongly stationary system is the direct product of a Bernoulli
system and a distal  strongly stationary system. By Theorem
\ref{T:C} a distal strongly stationary system coincides with its
$\C$-factor, so Theorem \ref{T:Bryna} shows that almost every
ergodic component of the distal factor of the system is a
pro-nilsystem. Finally, by Proposition \ref{P:jenvey} the system
has no rational eigenvalue different than $1$, so a nontrivial set
of ergodic components cannot share the same rational eigenvalue
provided that it is different than $1$. It follows that almost
every ergodic component of the system is totally ergodic, which
completes the proof.
\end{proof}

We will discuss the structure of the global system in more detail
in the next section. We will see that there exist distal strongly
systems with nonaffine ergodic components. This new class of
examples will allow a complete classification.

\subsection{New examples and structure theorem.}\label{SS:last}
Let $\sigma$ be a stationary measure on the sequence space $I^\Z$.
We define a new measure $\sigma_{av}$ on $I^\Z$ by averaging the
statistics of $\sigma$ along arithmetic progressions. More
precisely if $\F_0$ is the algebra of bounded $x_0$-measurable
functions and $S$ denotes the shift transformation, we define the
measure $\sigma_{av}$ on cylinder sets by
\begin{equation}\label{E:av}
\int f_0\ \bS f_1 \ \cdots \ \bS^k f_k \
d\sigma_{av}=\lim_{N\to\infty}\frac{1}{N} \sum_{n=1}^N \int f_0 \
\bS^nf_1\ \cdots \ \bS^{kn} f_k \ d\sigma,
\end{equation}
for $f_i\in\F_0$ (the limit exists by Theorem \ref{T:Bryna}). We
then extend $\sigma_{av}$ to the whole sequence space.

\begin{theorem}\label{T:averaging}
Let $\sigma$ be a totally  ergodic stationary measure. Then the
measure  $\sigma_{av}$ is   strongly stationary.
\end{theorem}
\begin{proof}
A direct computation proves stationarity. Since $S$ is totally
ergodic by Theorem \ref{C:Ziegler} we have
 $$
 \lim_{N\to\infty}  \frac{1}{N}
\sum_{n=1}^N \int f_0 \ \bS^{n}f_1\  \cdots \ \bS^{kn} f_k \
d\sigma= \lim_{N\to\infty}  \frac{1}{N} \sum_{n=1}^N \int f_0 \
\bS^{rn}f_1\ \cdots \ \bS^{rkn} f_k \ d\sigma
$$
for all $r,k \in\N$, $f_i\in\F_0$. Hence
$$
\int f_0\ \bS f_1 \ \cdots \ \bS^k f_k \ d\sigma_{av}= \int f_0\
\bS^r f_1 \ \cdots \ \bS^{rk} f_k \ d\sigma_{av}
$$
for all $r,k \in\N$, $f_i\in\F_0$. Since the subalgebra $\F_0$ is
$S$-generating the measure $\sigma_{av}$ is strongly stationary.
\end{proof}

 We remark that under the hypothesis of the previous theorem
the strongly stationary measure $\sigma_{av}$ can be shown to be
extremal. We will not use this fact so we omit its proof.

\begin{proposition}\label{P:av}
Let $\sigma$ be an extremal strongly stationary measure with
ergodic decomposition $\sigma=\int \sigma_t \ d\lambda(t)$. Then
for $\lambda$ almost every $t$ we have $\sigma=(\sigma_t)_{av}$.
\end{proposition}
\begin{proof}
Strong stationarity gives
\begin{align*}
\int f_0\ \bS f_1 \ \cdots \ \bS^k f_k \ d\sigma=&
\lim_{N\to\infty}
\frac{1}{N} \sum_{n=1}^N \int f_0 \ \bS^{n}f_1\  \cdots \ \bS^{kn} f_k \ d\sigma\\
=&\int \left( \lim_{N\to\infty}  \frac{1}{N} \sum_{n=1}^N \int f_0
\ \bS^{n}f_1\ \cdots \ \bS^{kn} f_k \ d\sigma_t \right)
d\lambda(t)
\end{align*}
for all $k \in\N$, $f_i\in\F_0$. It follows that
$$
\sigma=\int (\sigma_t)_{av}\ d\lambda(t).
$$
By Theorem \ref{P:jenvey} almost every ergodic component of $S$ is
totally ergodic so Theorem \ref{T:averaging} gives that  the
measures $(\sigma_t)_{av}$ are strongly stationary for
$\lambda$-a.e. $t$. Since $\sigma$ is extremal we  have that
$\sigma=(\sigma_t)_{av}$ for $\lambda$-a.e. $t$.
\end{proof}

We will use Theorem \ref{T:averaging} to construct an ample supply
of strongly stationary systems. We briefly describe the strategy.
Starting with an arbitrary invertible totally ergodic measure
preserving system ${\bf X}$  we first consider its sequence space
representation with respect to $\F=L^\infty(\bX)$
(Proposition~\ref{P:model2}). This representation is determined by
a stationary measure $\sigma$ on $I^\Z$, so ${\bf X}$ is
isomorphic to the system $(I^\Z,\B^\Z,\sigma,S)$ where $S$ is the
shift transformation. Let $\phi\colon X\to I^\Z$ be the
isomorphism ($\phi(\F)=\F_0=x_0$-measurable functions). We
construct the strongly stationary measure $\sigma_{av}$ as in
(\ref{E:av}), that is, we define the measure $\sigma_{av}$ on
cylinder sets by
\begin{equation}\label{E:av1}
\int f_0\ \bS f_1 \ \cdots \ \bS^k f_k \
d\sigma_{av}=\lim_{N\to\infty}\frac{1}{N} \sum_{n=1}^N \int f'_0 \
\bT^nf'_1\ \cdots \ \bT^{kn} f'_k \ d\mu,
\end{equation}
where $f_i\in \F_0$ and $f'_i=f_i\circ\phi\in\F$. Finally, we give
an explicit description of the statistics of $\sigma_{av}$. This
way we recover all the examples mentioned in Section
\ref{S:examples} and we also construct some new ones.

 \noindent {\bf Examples.} (i) Suppose that $\bX$ is a weak mixing
system. Using the multiple weak mixing theorem (\cite{Fu2}, page
86)  we can check that the resulting strongly stationary measure
defined by \eqref{E:av1} is a Bernoulli measure.

(ii) Suppose that $\bX$ is the system induced by  an ergodic
rotation on $\T$ with the Haar measure $m$. If we compute the
limit in (\ref{E:av1}) we find that
$$
\int f_0\ \bS f_1 \ \cdots \ \bS^k f_k \ d\sigma_{av}= \int_{\T^2}
f'_0(y)\
 f'_1(y+x) \ \cdots \  f'_k(y+kx) \ dm(y)\ dm(x).
$$
 We can  check that $\sigma_{av}$ determines the sequence
space representation of the strongly stationary system $T(x,y)=
(x,y+x)$ on $\T^2$ with the Haar measure $m$ (with respect to the
algebra generated by the exponentials in $y$).

(iii) Suppose that $\bX$ is the system induced by an affine
transformation on $\T^2$ with the Haar measure defined by
$T(x,y)=(x+a,y+x)$, where $a$ is irrational. We can check that the
resulting strongly stationary measure defined by \eqref{E:av1}
determines the sequence space representation of the system
$T'(x,y,z)=(x,y+x,z+y)$ on $\T^3$ with the Haar measure (with
respect to the algebra generated by the exponentials in $z$).

(iv) Suppose that $\bX$ is  an order $l$ totally ergodic nilsystem
defined on $X=G/\gG$. Using Theorem \ref{T:Ziegler} we see that
the resulting strongly stationary measure $\sigma_{av}$
constructed by \eqref{E:av1} is defined on cylinder sets as
follows
\begin{align}\label{E:nilsystem}
&\int f_0\ \bS f_1 \ \cdots \ \bS^k f_k \ d\sigma_{av}=
\int_{G/\gG}\int_{H_k/\Delta_k} f'_0(x\gG)f'_1(xy_1\gG)\ \cdots \
f'_k(xy_k\gG) \ d\nu_{H_k}  d\mu,
\end{align}
where  $H_k$, $\Delta_k$, $\nu_{H_k}$ are defined as in Theorem
\ref{T:Ziegler}. Then the system $(I^\Z,\B^\Z,\sigma_{av},S)$  is
strongly stationary with respect to $\F_0$.

(v) Suppose that $\bX$ is  a totally ergodic pro-nilsystem defined
on the inverse limit $X$ of the nilmanifolds $X_i$ ($X$ has to be
connected since it supports a totally ergodic pro-nilsystem). Then
the resulting strongly stationary measure $\sigma_{av}$
constructed by \eqref{E:av1} is defined on $\phi(X_i)$ by
\eqref{E:nilsystem}. Since $I^\Z=\bigcup_i \phi(X_i)$ this
uniquely determines the measure $\sigma_{av}$ on $I^\Z$. Observe
that the resulting strongly stationary system
$(I^\Z,\B^\Z,\sigma_{av},S)$ depends only on the pro-nilmanifold
$X$. We call it \emph{the strongly stationary system associated to
the pro-nilmanifold} $X$.

We remark that in example (iv) if $X$ is the Heisenberg
nilmanifold (see Section 6.1 example (i)) then the ergodic
components of the strongly stationary system obtained are
nilrotations on $X$, and hence nonaffine (\cite{Fu3}, page 52).



This new set of examples enables us to completely determine the
structure of the general strongly stationary system. This is the
context of our main theorem:

\begin{theorem}\label{T:main}
Every extremal strongly stationary system $(\bX,\F)$ is isomorphic
to the direct product of a Bernoulli system and a strongly
stationary system associated to some pro-nilmanifold.
\end{theorem}
\begin{proof}

Suppose that the ergodic components of $\bX$ are the systems
$\bX_t$. Let $\sigma$, $\sigma_t$ be the measures that determine
the sequence space representations of $\bX$ and $\bX_t$ with
respect to $L^\infty(\bX)$ and $L^\infty(\bX_t)$. Then
$\sigma=\int \sigma_t\ d\lambda$ is the ergodic decomposition of
$\sigma$. Since $\sigma$ is an extremal strongly stationary
measure by Proposition \ref{P:av}
  we have $\sigma=(\sigma_t)_{av}$ for $\lambda$-a.e. $t$. Theorem
 \ref{T:components} implies that such a measure $\sigma_t$ has the form
 $\rho \times \tau$
 where $\rho$ and  $\tau$  determine the sequence space representation
of a Bernoulli system  and a totally ergodic pro-nilsystem ${\bf
N}$ correspondingly. Then $\sigma=(\rho \times
 \tau)_{av}=\rho_{av} \times \tau_{av}$. The measure  $\rho_{av}$
 induces a Bernoulli system ${\bf B}$ (example (i))
and   the measure $\tau_{av}$  induces a strongly stationary
system ${\bf P}$ associated to the pro-nilmanifold $N$ (example
(v)). Hence, ${\bf X}$ is isomorphic to the direct product ${\bf
B}\times {\bf P}$ and the result follows.
\end{proof}

{\bf Acknowledgements.} Most of the results presented in this
paper are  contained in the author's PhD thesis. I wish to thank
my advisor Yitzhak Katznelson for his help and encouragement and
for many beneficial discussions. I would also like to thank Hillel
Furstenberg, Bernard Host,  and Bryna Kra for several helpful
suggestions.


\end{document}